\documentclass[12pt,twoside]{amsart}
\usepackage{amsthm,amsfonts,amssymb,amscd,euscript}

\usepackage[matrix,arrow]{xy}

\author{Florin Ambro} 
\address{RIMS, Kyoto University\\
Kyoto 606-8502, Japan.}
\email{ambro@kurims.kyoto-u.ac.jp}

\newcommand{\isoto}{{\overset{\sim}{\rightarrow}}}

\newcommand{\Q}{{\mathbb Q}}
\newcommand{\Z}{{\mathbb Z}}
\newcommand{\N}{{\mathbb N}}
\newcommand{\R}{{\mathbb R}}

\newcommand{\cL}{{\mathcal L}}

\newcommand{\cO}{{\mathcal O}}
\newcommand{\cR}{{\mathcal R}}

\newcommand{\bA}{{\mathbf A}}

\newcommand{\bD}{{\mathbf D}}

\newcommand{\bH}{{\mathbf H}}
\newcommand{\bK}{{\mathbf K}}

\newcommand{\bM}{{\mathbf M}}

\newcommand{\bP}{{\mathbf P}}

\newcommand{\Bsl}{\operatorname{Bsl}}
\newcommand{\codim}{\operatorname{codim}}
\newcommand{\BDiv}{\operatorname{BDiv}}
\newcommand{\CDiv}{\operatorname{CDiv}}

\newcommand{\LCS}{\operatorname{LCS}}

\newcommand{\mult}{\operatorname{mult}}

\newcommand{\rank}{\operatorname{rank}}

\newcommand{\Supp}{\operatorname{Supp}}

\newcommand{\WDiv}{\operatorname{WDiv}}

\theoremstyle{plain}
\newtheorem{thm}{Theorem}[section]

\newtheorem{lem}[thm]{Lemma}
\newtheorem{cor}[thm]{Corollary}

\newtheorem{prop}[thm]{Proposition}

\theoremstyle{definition}

\newtheorem{exmp}[thm]{Example}
\newtheorem{exmps}[thm]{Examples}

\theoremstyle{remark}

\newenvironment{sketch}{\begin{proof}[Sketch of proof]}{\end{proof}}

\begin{document}

\bibliographystyle{amsalpha+}
\title[A semiampleness criterion]
{A semiampleness criterion}
%\date{August 4, 2004}
\maketitle

\begin{abstract} 
We give a criterion for a real divisor to be rational
and semiample.
\end{abstract}

%\tableofcontents

%%%%%%%%%%%%%%%%%%%%%%%%%%%%%%%%%%%%%%%%%%%%%%%%%%%%%%%%%%%%%%%%%%%%%%
\setcounter{section}{-1}
%%%%%%%%%%%%%%%%%%%%%%%%%%%%%%%%%%%%%%%%%%%%%%%%%%%%%%%%%%%%%%%%%%%%%%
%%% Document name: semiacrt.tex
%%% Last modified: Sat Nov 20 00:20:32 JST 2004
%%%%%%%%%%%%%%%%%%%%%%%%%%%%%%%%%%%%%%%%%%%%%%%%%%%%%%%%%%%%%%%%%%%%%%

%%%%%%%%%%%%%%%%%%%%%%%%%%%%%%%%%%%%%%%
%%%%%%%%%%%%%%%%%%%%%%%%%%%%%%%%%%%%%%%

\section{Introduction}

%%%%%%%%%%%%%%%%%%%%%%%%%%%%%%%%%%%%%%%
%%%%%%%%%%%%%%%%%%%%%%%%%%%%%%%%%%%%%%%

\footnotetext[1]{This work was supported through a 
Twenty-First Century COE Kyoto Mathematics Fellowship.
}
\footnotetext[2]{1991 Mathematics Subject Classification. 
Primary: 14C20, 14E30. Secondary: 14N30.}

The aim of this paper is to establish the following criterion 
for the rationality and semiampleness of a nef divisor with 
real coefficients (we work over an algebraically closed field 
of characteristic zero):

\begin{thm}\label{M} Let $X$ be a proper nonsingular 
algebraic variety containing a simple normal crossings 
$\sum_l E_l$. 
Let $B=\sum_l b_l E_l$ and $D=\sum_l d_l E_l$, with 
$d_l,b_l\in \R$, be $\R$-divisors on $X$ satisfying 
the following properties: 
\begin{itemize}
\item[(1)] $D$ is nef and $rD-(K_X+B)$ is nef and big for 
some $r\in \Q$.
\item[(2)] $b_l<1$ for every $l$.
\item[(3)] For every integer $j\ge 1$, the following 
inclusion holds
$$
H^0(X,\sum_l \lceil jd_l-b_l\rceil E_l) \subseteq 
H^0(X,\sum_l \lfloor jd_l\rfloor E_l).
$$ 
For a real number $d$, $\lfloor d\rfloor$ and 
$\lceil d\rceil$ are the largest and smallest integers, 
respectively, such that $\lfloor d\rfloor\le d\le \lceil d\rceil$.
\end{itemize}
Then there exists a positive integer $m$ such that 
$md_l\in \Z$ for every $l$, and the linear 
system $\vert mD\vert$ is base point free.
\end{thm}

This criterion is also valid for singular varieties, being 
birational in nature (see Theorem~\ref{fr}).
It is a generalization of the Base Point Free Theorem 
(Kawamata-Shokurov~\cite{KMM}), and it also extends the known fact 
that the existence of the Cutkosky-Kawamata-Moriwaki decomposition 
for a log canonical divisor of general type implies the finite 
generation of the log canonical ring (Kawamata~\cite{lcd}).

Our interest in such a criterion lies in its connection to
the finite generation of Shokurov's FGA-algebras (\cite{PLflips}, 
Conjecture 4.39). Shokurov has recently reduced the existence of 
log flips in dimension $d+1$ to the Log Minimal Model Program in 
dimension $d$, and the finite generation of $d$-dimensional FGA-algebras 
(\cite{PLflips}, Corollary 1.5). Thus, provided that FGA-algebras 
are finitely generated, and this may be established even assuming 
the validity of the Log Minimal Model Program in the same dimension, 
this would give an inductive proof of the existence of log flips.
The finite generation of FGA-algebras is known in dimension one
and two (\cite{PLflips}, Main Theorem 1.7).

For the purpose of this introduction, an FGA-algebra is an 
$\N$-graded sub-algebra of a divisorial algebra
$$
\cL\subseteq \bigoplus_{i=0}^\infty H^0(X,iD),
$$
where $D$ is some divisor on a Fano variety $X$, such 
that $\cL$ is {\em asymptotically saturated with respect 
to} $X$. More precisely, the normalization of any such 
sub-algebra admits a representation
$$
\bar{\cL}=\bigoplus_{i=0}^\infty H^0(X_i,M_i),
$$
where each $M_i$ is a free Cartier divisor on some resolution 
of singularities $\mu_i\colon X_i\to X$ such that both 
divisors $K_{X_i}-\mu_i^*(K_X)$ and $M_i$ are supported 
by a simple normal crossings divisor on $X_i$. Then 
asymptotic saturation means that the following inclusions 
hold, for every $i,j\ge 1$:
$$
H^0(X_i,\lceil K_{X_i}-\mu_i^*(K_X)+j\frac{M_i}{i}\rceil)
\subseteq H^0(X_j,M_j).
$$ 
FGA-algebras are not well understood at the moment, the main
difficulty being in dealing with the a priori 
infinitely many models $X_i$. The successful two dimensional
approach to FGA-algebras is based on the Canonical Confinement 
Conjecture (\cite{PLflips}, Conjecture 6.14). Alternatively, 
based on the experience from the proof of the Base Point Free 
Theorem, Shokurov suggested that an inductive method may be 
possible in a larger class of algebras: $X$ may no longer be 
a Fano variety, but it should be related to $\cL$ by some kind 
of adjoint-type property (\cite{PLflips}, Remark 4.40.3). As 
an application of Theorem~\ref{M}, we give two such examples of 
adjoint-type properties which are sufficient for finite 
generation:

\begin{thm} Let $X$ be a proper nonsingular algebraic variety 
containing a simple normal crossings divisor $\sum_l E_l$.
Let $D=\sum_l d_l E_l$ be a nef a big $\R$-divisor, with 
associated divisorial algebra
$$
\cR_X(D)=\bigoplus_{i=0}^\infty H^0(X,iD).
$$
Assume that there exists an $\R$-divisor $B=\sum_l b_l E_l$,
with $b_l<1$ for every $l$, such that 
\begin{itemize}
\item[(1)] $\cR_X(D)$ is asymptotically saturated with respect 
to $(X,B)$.
\item[(2)] The $\R$-divisor $rD-(K_X+B)$ is nef, for some $r\in \Q$. 
\end{itemize}
Then $\cR_X(D)$ is a finitely generated algebra.
\end{thm}

\begin{thm}\label{onemodel} 
Let $X$ be proper nonsingular algebraic variety containing
a simple normal crossings divisor $\sum_l E_l$. Let 
$M_i=\sum_l m_{i,l}E_l$ be divisors on $X$, for $i\ge 0$, 
satisfying the following properties:
\begin{itemize}
\item[(i)] The linear system $\vert M_i\vert$ is 
base point free, for every $i$.
\item[(ii)] $M_0=0$, $M_i+M_j\le M_{i+j}$ for every $i,j$.
\item[(iii)] The limit $d_l=\lim_{i\to \infty} \frac{1}{i}m_{i,l}$ 
exists, for every $l$. Set $D=\sum_l d_l E_l$.
\end{itemize}
The sequence $(M_i)_{i\ge 0}$ defines an algebra 
$$
\cL=\bigoplus_{i=0}^\infty H^0(X,M_i)\subseteq \cR_X(D).
$$ 
Assume moreover that there exists an $\R$-divisor 
$B=\sum_l b_l E_l$, with $b_l<1$ for every $l$, such that
\begin{itemize}
\item[(1)] $\cL$ is asymptotically saturated with respect
to $(X,B)$. Equivalently,
$H^0(X,\sum_l \lceil j\frac{m_{i,l}}{i}-
b_l\rceil E_l) \subseteq H^0(X,M_j)$ for every $i,j\ge 1$.
\item[(2)] The $\R$-divisor $rD-(K_X+B)$ is nef and big,
for some $r\in \Q$.
\end{itemize}
Then $\cL$ is a finitely generated algebra.
\end{thm}

The paper is organized as follows. In Section 1, we collect
some basic definitions and properties on log pairs, b-divisors 
and functional algebras. In Section 2, we establish the 
semiampleness criterion. The proof is quite similar to that 
of the Base Point Free Theorem, and we use a nonvanishing 
result for real divisors due to Kawamata (\cite{lcd}, Theorem 3). 
In Section 3, we apply this criterion to the finite generation of 
some functional algebras.

%\clearpage
%%%%%%%%%%%%%%%%%%%%%%%%%%%%%%%%%%%%%%%
%%%%%%%%%%%%%%%%%%%%%%%%%%%%%%%%%%%%%%%

\section{Preliminary}

%%%%%%%%%%%%%%%%%%%%%%%%%%%%%%%%%%%%%%%
%%%%%%%%%%%%%%%%%%%%%%%%%%%%%%%%%%%%%%%

A {\em variety} is a reduced and irreducible 
separable scheme of finite type, defined over an 
algebraically closed field $k$ of characteristic zero. 
A {\em contraction} is a proper morphism of normal varieties
$f\colon X \to Y$ with $\cO_Y=f_*\cO_X$.

%%%%%%%%%%%%%%%%%%%%%%%%%%%%%%%%%%
%%%%%%%%%%%%%%%%%%%%%%%%%%%%%%%%%%

\subsection{Divisors}

%%%%%%%%%%%%%%%%%%%%%%%%%%%%%%%%%%%
%%%%%%%%%%%%%%%%%%%%%%%%%%%%%%%%%%%
Let $X$ a normal variety. A Weil divisor on $X$ is an
element of the free abelian group $\WDiv(X)$ generated 
by codimension one cycles on $X$. A non-zero rational
function $a\in k(X)^\times$ defines the principal Weil
divisor 
$$
(a)=\sum_{\codim(P,X)=1}\mult_P(a)P.
$$
Let $L\in \{\Z,\Q,\R\}$. An $L$-Weil divisor is an element 
of 
$$
\WDiv(X)_L:=\WDiv(X)\otimes_\Z L. 
$$
We say that $D_1, D_2\in \WDiv(X)_\R$ are 
$L$-{\em linearly equivalent}, denoted $D_1\sim_L D_2$, 
if there exist finitely many $q_i\in L$ and $a_i\in k(X)^\times$ 
such that $D_1-D_2=\sum_i q_i(a_i)$. 
We say that $D\in \WDiv(X)_L$ is {\em $L$-Cartier} if $D\sim_L 0$ 
in a neighborhood of each point of $X$. We denote by
$$
\CDiv(X)_L\subset \WDiv(X)_L
$$ 
the subgroup of $L$-Cartier $L$-Weil divisors. 
To each $D\in \WDiv(X)_\R$ we can associate a coherent 
$\cO_X$-sub-module $\cO_X(D)$ of the constant sheaf of 
rational functions $k(X)$ on $X$, whose sections on 
open subsets $U$ of $X$ are given by
$$
H^0(U,\cO_X(D))=\{a\in k(X)^\times; ((a)+D)\vert_U\ge 0\}\cup \{0\}.
$$
For $D\in \WDiv(X)_\R$, its {\em round down} $\lfloor D\rfloor$,
{\em round up} $\lceil D\rceil$, {\em fractional part} 
$\{D\}$ are defined componentwise (for $d\in \R$, $\lfloor d\rfloor$
and $\lceil d\rceil$ are the unique integers in the half-open intervals 
$(d-1,d]$ and $[d,d+1)$, respectively, and $\{d\}=d-\lfloor d\rfloor$). 

Let $\pi\colon X\to S$ be a proper surjective morphism. 
We say that $D\in \WDiv(X)_\R$ is 
\begin{itemize}
\item[(i)] {\em free}/$S$ if $D$ is Cartier and 
the natural map $\pi^*\pi_*\cO_X(D)\to \cO_X(D)$ is 
surjective.
\item[(ii)] {\em nef}/$S$ if $D$ is $\R$-Cartier and 
$D\cdot C\ge 0$ for every proper curve $C\subset X$ such 
that $\pi(C)$ is a point.
\item[(iii)] {\em ample}/$S$ if $\pi$ is projective and 
the numerical class of $D$ belongs to the real cone 
generated by the numerical classes of $\pi$-ample Cartier 
divisors.
\item[(iv)] {\em semiample}/$S$ if there exists a 
contraction $f\colon X\to Y$, defined over $S$, and an
ample/$S$ $\R$-divisor $H$ on $Y$ such that $D\sim_\R f^*H$. 
If $D$ is rational, this is equivalent to the existence
of a positive integer $m$ such that $mD$ is free/$S$.
\item[(v)] {\em big}/$S$ if there exists $C>0$ such that
$\rank \pi_*\cO_X(mD)\ge Cm^d$ for $m$ sufficiently 
large and divisible, where $d$ is the dimension of the 
general fiber of $\pi$. 
\end{itemize}

%%%%%%%%%%%%%%%%%%%%%%%%%%%%%%%%%%
%%%%%%%%%%%%%%%%%%%%%%%%%%%%%%%%%%

\subsection{B-divisors}

%%%%%%%%%%%%%%%%%%%%%%%%%%%%%%%%%%%
%%%%%%%%%%%%%%%%%%%%%%%%%%%%%%%%%%%
(V.V. Shokurov~\cite{Logmodels}). Let $X$ be a normal
variety. The group of {\em b-divisors} of $X$ is defined as
the projective limit
$$
\BDiv(X)=\projlim_{Y\to X} \WDiv(Y),
$$ 
where the index set is the set of isomorphism classes of 
birational contractions $Y\to X$, with the following order: 
$Y_1\ge Y_2$ if the induced proper
birational map $g\colon Y_1\dashrightarrow Y_2$ is a morphism.
In this case, we have a natural push-forward homomorphism
$g_*\colon \WDiv(Y_1)\to \WDiv(Y_2)$.
Note that given two birational contractions $Y_1\to X$ and 
$Y_2\to X$, the normalization of the graph of the induced 
birational map $Y_1\dashrightarrow Y_2$ defines a 
third birational contraction $Y_3\to X$ with $Y_3\ge Y_1$ 
and $Y_3\ge Y_2$.

A {\em prime b-divisor} $E$ of $X$ is obtained as follows:
let $\mu\colon Y\to X$ be a birational contraction, and let 
$E$ be a prime divisor on $Y$. For another birational contraction
$Y'\to X$, let $g\colon Y \dashrightarrow Y'$ be the induced 
rational map. If $g$ induces an isomorphism between the generic
point of $E$ and the generic point of a prime divisor $E'$ on $Y'$,
we set $\bD_{Y'}=E'$. Otherwise, set $\bD_{Y'}=0$. The 
family $\{\bD_{Y'}\}$ is a prime b-divisor of $X$, denoted
by $E$. The center of $E$ on $X$, denoted by $c_X(E)$, is
the closed subset $\mu(E)$ of $X$.
It is clear that prime b-divisors of $X$ are in one-to-one 
correspondence with geometric (discrete) valuations of $k(X)$ 
which have non-empty center on $X$. In this terminology, a 
b-divisor of $X$ is a $\Z$-valued function on the set of all
prime b-divisors of $X$, denoted as a formal sum
$$
\bD=\sum_E \mult_E(\bD) E,
$$
such that the sum has finite support on some (hence any)
birational model of $X$. 

A proper birational map $\Phi\colon X_1\dashrightarrow X_2$ 
(i.e., the normalization of the graph of $\Phi$ induces 
birational contractions to $X_1$ and $X_2$), induces a 
push-forward isomorphism
$
\Phi_*\colon \BDiv(X_1)\isoto \BDiv(X_2).
$

For $L\in \{\Z,\Q,\R\}$, an {\em $L$-b-divisors} of $X$ is 
an element of 
$$
\BDiv(X)_L:=\BDiv(X)\otimes_\Z L.
$$

\begin{exmps}
\begin{itemize}
\item[(1)] A rational function $a \in k(X)^\times$ defines the
following b-divisor of $X$ 
$$
\overline{(a)}=\{(f^*a)\}_{f\colon Y\to X}.
$$
\item[(2)] Let $Y\to X$ be a birational contraction and 
let $D\in \CDiv(Y)_\R$. The family 
$\overline{D}=\{f^*D\}_{f\colon Y'\to Y}$ is an 
$\R$-b-divisor of $X$, called the {\em Cartier closure}
of $D$.
\item[(3)] Let $\omega\in \Omega^{\dim(X)}_{k(X)/k}$ be a 
top rational differential form of $X$. The associated family 
of Weil divisors $\bK=\{(f^*\omega)\}_{f\colon Y\to X}$ is called 
the {\em canonical b-divisor} of $X$. It is uniquely defined 
up to linear equivalence.
\end{itemize}
\end{exmps}
Let $\bD$ be an $\R$-b-divisor of $X$. The round up 
$\lceil \bD\rceil\in \BDiv(X)$ is defined componentwise.
The restriction of $\bD$ to an open subset $U\subset X$ is 
a well defined $\R$-b-divisor of $U$, denoted by $\bD\vert_U$.
Also, $\bD$ defines an $\cO_X$-module $\cO_X(\bD)$ whose sections 
on an open subset $U\subseteq X$ are given by
$$
H^0(U,\cO_X(\bD))=\{a\in k(X)^\times; 
(\overline{(a)}+\bD)\vert_U\ge 0\}\cup\{0\}.
$$
An $\R$-b-divisor $\bD$ of $X$ is called {\em $L$-b-Cartier} 
if there exists a birational contraction $f\colon X'\to X$ 
such that $\bD_{X'}$ is $L$-Cartier and $\bD=\overline{\bD_{X'}}$.
Note that $\cO_X(\bD)=f_*\cO_{X'}(\bD_{X'})$ in this case,
in particular $\cO_X(\bD)$ is a coherent $\cO_X$-module.

Let $\pi\colon X\to S$ be a proper morphism. An $\R$-b-divisor 
$\bD$ of $X$ is called {\em b-free}/$S$ ({\em b-nef}/$S$, 
{\em b-semiample}/$S$, {\em b-nef and big}/$S$) if there exists 
a birational contraction $X'\to X$ such that $\bD=\overline{\bD_{X'}}$ 
and $\bD_{X'}$ is free/$S$ (nef/$S$, semiample/$S$, nef/$S$ and 
big/$S$), where $X'\to S$ is the induced proper morphism.

\begin{lem}[\cite{PLflips}, Proposition 4.15]\label{mobiledef} 
Let $\pi\colon X\to S$ be a contraction. 
Let $k(X)$ be the constant sheaf of rational functions on $X$ and
let $\cL\subset \pi_*k(X)$ be a non-zero coherent $\cO_X$-module.
Then there exists a unique b-free/$S$ b-divisor $\bM$ of $X$ 
such that 
\begin{itemize}
\item[(i)] $\cL\subseteq \pi_*\cO_X(\bM)$.
\item[(ii)] If $\bD\in \BDiv(X)_\R$ and
$\cL\subseteq \pi_*\cO_X(\bD)$, then $\bM\le \bD$.
\end{itemize}
We say that $\bM$ is the {\em mobile b-divisor} of $\cL$.
\end{lem}

\begin{sketch} Assume first that $S$ is affine, so that $\cL$ is
generated by finitely many global sections $a_1,\ldots,a_p$. The 
b-divisor $\bM$ is defined by
$$
\bM=\sup\{-\overline{(a)}; a\in H^0(S,\cL)\setminus 0\}
=\max_{i=1}^p (-\overline{(a_i)})
$$
In general, we may cover $S$ with finitely many affine open 
subsets $S_\alpha$, and the b-divisors $\bM_\alpha$
of $\pi^{-1}(S_\alpha)$ constructed above glue to a b-divisor 
$\bM\in \BDiv(X)$, which clearly satisfies (i) and (ii).

To show that $\bM$ is b-free/$S$, assume again that $S$ is
affine. By the definition of $\bM_X$, the linear system
$$
\{(a)+\bM_X;a\in H^0(S,\cL)\setminus 0\}
$$
has no fixed components. By Hironaka, there 
exists a birational contraction $\mu\colon Y\to X$ such that 
the induced linear system 
$$
\{(\mu^*a)+\bM_Y;a\in H^0(S,\cL)\setminus 0\}
$$
has no fixed points. Then $\bM_Y$ is free/$S$ and 
$\bM=\overline{\bM_Y}$.
\end{sketch}

\begin{lem}\label{neg} Let $\pi\colon X\to S$ be a contraction. 
Let $\bD\in \BDiv(X)_\R$ such that $\bD_X\le 0$ and $\bD$ is 
b-nef/$S$. Then:
\begin{itemize}
\item[(1)] $\bD\le 0$.
\item[(2)] There exists a closed subset $S^0\subsetneqq S$
such that for every prime b-divisor $E$ of $X$,
$\mult_E(\bD)\ne 0$ if and only if $\pi(c_X(E))\subseteq S^0$.
\end{itemize}
\end{lem}

\begin{proof} (1) Let $Y\to X$ be a birational contraction
such that $\bD=\overline{\bD_Y}$. Since $\bD_Y$ is nef/$S$, 
it is also nef/$X$.
By the Negativity Lemma~\cite[1.1]{3flips}, we obtain
$\bD_Y\le 0$. Therefore $\bD=\overline{\bD_Y}\le 0$.

(2) Let $Y\to X$ be a birational contraction such that 
$\bD=\overline{\bD_Y}$, and let $\nu\colon Y\to S$ be
the induced contraction. 
The restriction of the $\R$-divisor $-\bD_Y$ to the
general fibre of $\nu$ is zero, being anti-nef and effective.
Therefore $$S^0=\nu(\Supp(\bD_Y))$$ is a closed subset
of $S$ and $S^0\ne S$. Property (2) is equivalent
to the following equality of sets
$$
\Supp(\bD_Y)=\nu^{-1}(S^0).
$$
The inclusion $\Supp(\bD_Y)\subseteq \nu^{-1}(S^0)$
is clear. For the converse, fix a point $P\in S^0$
and assume by contradiction that $\nu^{-1}(P)$
is not included in $\Supp(\bD_Y)$. We may replace 
$Y$ by a resolution of singularities such that
$\nu^{-1}(P)$ is the union of finitely many
prime Cartier divisors $E_1,\ldots,E_n$ on $Y$.
By assumption, we have $\mult_{E_1}(\bD_Y)=0$. The fiber
$\nu^{-1}(P)=\bigcup_{i=1}^n E_i$ is connected, hence
after a possible reordering, we may assume that $E_j$ 
intersects $E_1\cup\cdots \cup E_{j-1}$, 
for every $2\le j\le n$.

Assume by induction that $\mult_{E_i}(\bD_Y)=0$ for
$i\le j-1$. There exists a proper curve 
$C\subset \bigcup_{i=1}^{j-1}E_i$ such that 
$C\cap E_j\ne \emptyset$ and $C\nsubseteq 
\Supp(\bD_Y)$. Since $-\bD_Y$ is effective and 
$\bD_Y\cdot C\ge 0$, we infer $\mult_{E_j}(\bD_Y)=0$.
By induction, we obtain $\mult_{E_j}(\bD_Y)=0$ for
$1\le j\le n$. However, $P\in S^0$ implies that 
$\Supp(\bD_Y)\cap \nu^{-1}(P)\ne \emptyset$.
In particular, there exists a proper curve 
$C\subseteq \nu^{-1}(P)$ such that $C\not\subseteq \Supp(\bD_Y)$
and $C\cap \Supp(\bD_Y)\ne \emptyset$. Therefore 
$\bD_Y\cdot C<0$. On the other hand, $\bD_Y$ is nef/$S$, 
hence $\bD_Y\cdot C\ge 0$. Contradiction.
\end{proof}

\begin{lem}\label{neg'} Let $X\to S$ be a contraction, and 
let $\bD\le\bH$ be $\R$-b-divisors of $X$ such that $\bD$ is 
b-nef/$S$ and $\bH$ is $\R$-b-semiample/$S$, that is there 
exists rational contraction over $S$
\[ \xymatrix{
X\ar[dr] \ar@{.>}[rr]^h &  & T \ar[dl] \\
       &      S               &  
} \]
and an ample/$S$ $\R$-divisor $A$ on $T$ such that
$\bH\sim_\R \overline{h^*(A)}$.
If $\bD=\bH$ over a neighborhood of a closed subset
$X_0\subset X$, then $\bD=\bH$ over a neighborhood of 
$h(X_0)\subset T$.
\end{lem}

\begin{proof} We may replace $X$ by a higher birational model
so that $h$ becomes a morphism. By assumption, $\bD-\bH\le 0$ is 
b-nef/$T$. By Lemma~\ref{neg}, there exists a closed subset
$T^0\subset T$ such that $\mult_E(\bD-\bH)\ne 0$ if
and only if $c_X(E)\subseteq h^{-1}(T^0)$, for every prime
b-divisor $E$ of $X$. By assumption, there exists an open
set $U\supset X_0$ such that $\mult_E(\bD-\bH)= 0$ for
every prime b-divisor $E$ of $X$ with $c_X(E)\cap U\ne 
\emptyset$. Therefore $U\subseteq X\setminus h^{-1}(T^0)$.
In particular, $X_0\cap h^{-1}(T^0)=\emptyset$. 
Therefore $V=T\setminus T^0$ is an open neighborhood of $h(X_0)$,
and 
$
\bD\vert_{h^{-1}(V)}=\bH\vert_{h^{-1}(V)}.
$ 
\end{proof}

%%%%%%%%%%%%%%%%%%%%%%%%%%%%%%%%%%
%%%%%%%%%%%%%%%%%%%%%%%%%%%%%%%%%%

\subsection{Log pairs}

%%%%%%%%%%%%%%%%%%%%%%%%%%%%%%%%%%%
%%%%%%%%%%%%%%%%%%%%%%%%%%%%%%%%%%%
A {\em log pair} $(X,B)$ is a normal variety $X$ 
endowed with an $\R$-Weil divisor $B$ such that $K+B$ 
is $\R$-Cartier. A {\em log variety} is a log pair 
$(X,B)$ such that $B$ is effective.

Let $(X,B)$ be a log pair and let $\mu\colon Y\to X$ be
a birational contraction. Choose a top rational differential 
form $\omega\in \Omega^{\dim(X)}_{k(X)/k}$ with $K=(\omega)$. 
The Weil divisor $K_Y=(\mu^*\omega)$ is a canonical divisor 
of $Y$, and there exists a unique $\R$-Weil divisor $B_Y$
on $Y$ satisfying the adjunction formula
$$
\mu^*(K+B)=K_Y+B_Y.
$$
We say that $(Y,B_Y)$ is the {\em crepant log pair structure}
induced by $(X,B)$ via $\mu$. It is independent of the choice
of $\omega$ and in fact it is independent of the choice of
$K$ in its linear equivalence class. There exists a unique
$\R$-b-divisor $\bA(X,B)$ of $X$, called the 
{\em discrepancy $\R$-b-divisor} of the log pair $(X,B)$,
whose trace on birational modifications $Y\to X$ is 
$$
\bA(X,B)_Y=-B_Y.
$$
A log pair $(X,B)$ is said to have {\em Kawamata log terminal 
singularities}
if $\mult_E(\bA(X,B))>-1$ for every geometric valuation
$E$ of $X$. We denote by $\LCS(X,B)$ the largest closed
subset of $X$ on whose complement $(X,B)$ has Kawamata
log terminal singularities. 

For a log pair $(X,B)$, the $\cO_X$-module 
$\cO_X(\lceil \bA(X,B)\rceil)$ is coherent. Indeed,
let $\mu\colon Y\to X$ be a resolution of singularities such 
that $B_Y=\sum_l b_l E_l$ and $\sum_l E_l$ is a simple normal 
crossings divisor. Then (see the proof of~\cite{PLflips}, 
Proposition 4.46)
$$
\cO_X(\lceil \bA(X,B)\rceil)=\mu_*\cO_Y(\lceil -B_Y\rceil).
$$
Also, $\LCS(X,B)=\mu(\LCS(Y,B_Y))$ and 
$\LCS(Y,B_Y)=\bigcup_{b_l\ge 1}E_l$.

%%%%%%%%%%%%%%%%%%%%%%%%%%%%%%%%%%
%%%%%%%%%%%%%%%%%%%%%%%%%%%%%%%%%%

\subsection{Functional algebras}

%%%%%%%%%%%%%%%%%%%%%%%%%%%%%%%%%%%
%%%%%%%%%%%%%%%%%%%%%%%%%%%%%%%%%%%
(V.V. Shokurov~\cite{PLflips})
Let $\pi\colon X\to S$ be a contraction.
A {\em functional algebra of} $X/S$ is an $\N$-graded 
$\cO_S$-subalgebra 
$$
\cL \subset \bigoplus_{i=0}^\infty \pi_*k(X)
$$ 
such that $\cL_0=\cO_S$ and $\cL_i$ is a coherent 
$\cO_S$-module for every $i\ge 1$. Here, $k(X)$ is the constant 
sheaf of rational functions on $X$, and the multiplication in 
$\bigoplus_{i=0}^\infty \pi_*k(X)$ is induced by that of $k(X)$.

Denote by $\N(\cL)$ the semigroup of positive integers $i$ 
with $\cL_i\ne 0$. For $i\in \N(\cL)$, let 
$\bM_i$ be the mobile b-divisor of $\cL_i$ (see Lemma~\ref{mobiledef})
and set $\bD_i=\bM_i/i$. By construction, $\bM_i$ is b-free/$S$ 
and $\bD_i$ is a b-semiample/$S$ $\Q$-b-divisor. The sequences 
$(\bM_i)_{i\in \N(\cL)}$ and $(\bD_i)_{i\in \N(\cL)}$ are called 
the {\em mobile sequence} and {\em characteristic sequence} of $\cL$.
Since $\cL$ is an algebra, the following properties hold:
\begin{itemize}
\item[(i)] $\bM_i+\bM_j\le \bM_{i+j}$.
\item[(ii)] $\bM_0=0$.
\end{itemize}
In particular, $\bigoplus_{i\in \{0\}\cup\N(\cL)} \pi_*\cO_X(\bM_i)$ 
is a functional algebra of $X/S$, which is the normalization of $\cL$
(\cite{PLflips}, Proposition 4.15).

The functional algebra $\cL$, is called {\em bounded} if the
supremum 
$$
\bD=\sup_{i\in \N(\cL)}\bD_i
$$
is an $\R$-b-divisor of $X$, called the {\em characteristic limit} 
of $\cL$. In this case, the convexity property (i) implies that
the supremum is a limit
$$
\bD=\lim_{\N(\cL)\ni i\to\infty}\bD_i.
$$
The $\cO_S$-algebra $\cL$ is finitely generated if and only if 
$\bD_i=\bD$ for some $i\in \N(\cL)$ (\cite{PLflips}, Theorem 4.28).

Let $\bA\in \BDiv(X)_\R$. The functional algebra $\cL$ is
{\em asymptotically $\bA$-saturated} if there exists a positive 
integer $I$ such that $\cL_i\ne 0$ for $I\vert i$, and the 
following inclusion holds for $I\vert i,j$:
$$
\pi_*\cO_X(\lceil \bA+j\bD_i\rceil)\subseteq \pi_*\cO_X(\bM_j).
$$
We say that $\cL$ is {\em asymptotically saturated with respect
to a log pair} $(X,B)$ if $\cL$ is asymptotically 
$\bA(X,B)$-saturated.

\begin{exmp} Let $\pi\colon X\to S$ be a contraction.
\begin{itemize}
\item[(1)] Let $\bD$ be an $\R$-b-divisor of $X$. Then
$$
\cR_{X/S}(\bD)=\bigoplus_{i=0}^\infty \pi_*\cO_X(i\bD)
$$
is a bounded functional algebra of $X/S$.
\item[(2)] Let $(X,B)$ be a log variety structure on $X$ and let
$D\in \CDiv(X)_\Q$. Then the normal bounded functional algebra of 
$X/S$
$$
\cR_{X/S}(D)=\bigoplus_{i=0}^\infty \pi_*\cO_X(iD)
$$
is asymptotically $\bA(X,B)$-saturated.
\end{itemize}
\end{exmp}

%%%%%%%%%%%%%%%%%%%%%%%%%%%%%%%%%%%%%%%
%%%%%%%%%%%%%%%%%%%%%%%%%%%%%%%%%%%%%%%

\subsection{Kawamata's Nonvanishing}

%%%%%%%%%%%%%%%%%%%%%%%%%%%%%%%%%%%%%%%
%%%%%%%%%%%%%%%%%%%%%%%%%%%%%%%%%%%%%%%

For a real number $r$, the {\em round off} $\langle r\rangle$ is 
defined as the unique integer lying in the unit interval 
$[r-1/2,r+1/2)$. The following properties are easy to check:
\begin{itemize}
\item[(1)] $\lfloor r \rfloor \le \langle r\rangle 
\le \lceil r \rceil$. 
\item[(2)] $\langle -r\rangle=-\langle r\rangle$
unless $\{r\}=\frac{1}{2}$.
\item[(3)] $\langle r\rangle=0$ for $r\in 
(-\frac{1}{2},\frac{1}{2}]$.
\end{itemize}

\begin{lem}\cite{lcd}\label{app} 
Let $x\in \R^d$, let $I$ be a positive integer and
let $\epsilon>0$. Then there exists a positive multiple
$m$ of $I$ such that $mx_1\le\langle mx_1\rangle$ and
$
\max_{i=1}^d|mx_i-\langle mx_i\rangle|<\epsilon.
$
\end{lem}

For an $\R$-divisor $D$ on a normal variety $X$, the
{\em round off} $\langle D\rangle$ is defined componentwise.
Also, the {\em absolute value} of $D$ is defined by 
$$
\Vert D\Vert=\max\{\vert\mult_P(D)\vert; P\subset X,
\codim(P,X)=1\}.
$$

\begin{thm}\cite{lcd} \label{Knv} Let $X$ be a nonsingular 
variety and let $\pi\colon X\to S$ be a proper morphism. Let 
$D,B$ be $\R$-divisors on $X$, satisfying the following 
properties:
\begin{itemize}
\item[(1)] $D$ is nef/$S$.
\item[(2)] The $\R$-divisor $rD-(K+B)$ is nef/$S$ and big/$S$, 
for some $r\in \Q$.
\item[(3)] $\Supp\{B\}$ is a simple normal crossings divisor 
and $\lfloor B\rfloor\le 0$.
\end{itemize}
Then there exists $t_0,\epsilon>0$ such that
$
\pi_*\cO_X(\lceil -B\rceil + \langle tD\rangle)\ne 0
$
for every $t\in \R$ with $t\ge t_0$ and 
$\Vert \langle tD\rangle-tD\Vert <\epsilon$.
\end{thm}

\begin{cor} \label{Snv}
Let $(X,B)$ be a log pair with Kawamata log terminal
singularities and let $\pi\colon X\to S$ be a proper
morphism. Let $\bD$ be a b-nef/$S$ $\R$-b-divisor of 
$X$ such that the $\R$-b-divisor
$
r\bD-\overline{K+B}
$
is b-nef/$S$ and b-big/$S$ for some $r\in \R$. 

Then there exists $t_0, \epsilon>0$ and a birational 
contraction $Y\to X$ such that
$
\pi_*\cO_X(\lceil \bA(X,B)+t\bD \rceil) \ne 0
$
for every $t\in \R$ with $t\ge t_0$ and 
$\Vert t\bD_Y-\langle t\bD_Y\rangle \Vert <\epsilon$.
\end{cor}

\begin{proof} Let $Y\to X$ be a resolution of singularities
such that $\bD=\overline{\bD_Y}$, and $\Supp(B_Y)\cup\Supp(\bD_Y)$ 
is included in a simple normal crossings divisor 
$\sum_l E_l$ on $Y$. Let $\pi'\colon Y\to S$ be he induced
proper morphism. In particular,
$$
\pi_*\cO_X(\lceil \bA(X,B)+t\bD\rceil)=
\pi'_*\cO_Y(\lceil -B_Y+t\bD_Y\rceil)\mbox{ for } t\in \R.
$$
By Theorem~\ref{Knv}, there exist $t_0, \epsilon'>0$
such that 
$$
\pi'_*\cO_Y(\lceil -B_Y\rceil +\langle t\bD_Y\rangle) \ne 0
$$
for $t\ge t_0$ and $\Vert t\bD_Y-\langle t\bD_Y\rangle\Vert
<\epsilon'_0$. Choose $\epsilon\in (0,\epsilon')$ and
$$
\epsilon\le 
\min_l (1-\{\mult_{E_l}(B_Y)\}).
$$
Then for every $t\in \R$ with 
$
\Vert t\bD_Y-\langle t\bD_Y\rangle\Vert<\epsilon,
$
we have
$$
\lceil -B_Y\rceil +\langle t\bD_Y\rangle\le
\lceil -B_Y+t\bD_Y\rceil.
$$
Therefore 
$
\pi'_*\cO_Y(\lceil -B_Y+t\bD_Y \rceil)\ne 0
$
for $t\ge t_0$ and $\Vert t\bD_Y-\langle t\bD_Y\rangle \Vert
<\epsilon$.
\end{proof}

%%%%%%%%%%%%%%%%%%%%%%%%%%%%%%%%%%%%%%%
%%%%%%%%%%%%%%%%%%%%%%%%%%%%%%%%%%%%%%%

\section{A semiampleness criterion}

%%%%%%%%%%%%%%%%%%%%%%%%%%%%%%%%%%%%%%%
%%%%%%%%%%%%%%%%%%%%%%%%%%%%%%%%%%%%%%%

\begin{thm}\label{fr} Let $(X,B)$ be a log pair, let
$\pi\colon X\to S$ be a proper morphism and let $D$ 
be a nef/$S$ $\R$-divisor on $X$, satisfying the 
following properties:
\begin{itemize}
\item[(i)] $rD-(K+B)$ is nef/$S$ and big/$S$, for 
some $r\in \Q$.
\item[(ii)] $\pi_*\cO_X(\lceil \bA(X,B)+j\overline{D}\rceil)
\subseteq \pi_*\cO_X(jD)$ for $I|j$, where $I$ is some positive
integer.
\item[(iii)] For every point $P\in \LCS(X,B)$, there 
exists a non-zero rational function $a\in k(X)^\times$ 
and there exists $q\in \Q$ such that the $\R$-divisor 
$$
q(a)+D
$$ 
is effective in a neighborhood of $\pi^{-1}(\pi(P))$ and it
does not contain $P$ in its support.
\end{itemize}
Then $D$ is a semiample/$S$ $\Q$-divisor.
\end{thm}

\begin{proof} 
{\em Step 1}. There exists a positive multiple $I_1$ of $I$ 
such that 
$$
\pi_*\cO_X(mD)\ne 0 \mbox{ for } I_1\vert m.
$$
Indeed, if $(X,B)$ does not have Kawamata log terminal 
singularities, the claim follows from (iii). If $(X,B)$ has 
Kawamata log terminal singularities, we infer by Corollary~\ref{Snv} 
and Lemma~\ref{app} that there exists a multiple
$I_1$ of $I$ such that 
$$
\pi_*\cO_X(\lceil \bA(X,B)+I_1\overline{D}\rceil)\ne 0.
$$
By (ii), we obtain $\pi_*\cO_X(I_1D)\ne 0$. The same holds
for multiples of $I_1$.

{\em Step 2}. Let $(\bD_i)_i$ be the characteristic 
sequence of the functional algebra of $X/S$
$$
\cR_{X/S}(D)=\bigoplus_{i=0}^\infty \pi_*\cO_X(iD).
$$ 
By (iii), there exists a 
multiple $I_2$ of $I_1$ such that $\overline{D}-\bD_m$
is zero above a neighborhood of $\LCS(X,B)$, for every
$m\in I_2\Z_{>0}$. We replace $D,r$ by $I_2D, rI_2^{-1}$,
respectively, so that we can assume from now that 
$I=I_1=I_2=1$. We will also shrink $S$ to an affine
neighborhood of some fixed point of $S$, without further 
notice.

{\em Step 3}. Fix a positive integer $m$ such that 
$\bD_m\ne \overline{D}$. 
There exists a resolution $\mu\colon Y\to X$ with the 
following properties (cf.~\cite{KMM}, Corollary 0-3-6):
\begin{itemize}
\item[(1)] $\bD_m=\overline{(\bD_m)_Y}$. 
Denote $F=\mu^*D-\bD_{m,Y}$.
\item[(2)] $\Supp(B_Y)\cup \Supp(F)$ is included in a 
simple normal crossings divisor $\sum_\alpha G_\alpha$ 
on $Y$, where $K_Y+B_Y=\mu^*(K+B)$.
\item[(3)] 
$r\mu^*D-K_Y-B_Y-\sum_\alpha \epsilon_\alpha G_\alpha$
is an ample/$S$ $\R$-divisor, with $0<\epsilon_\alpha\ll 1$.
\end{itemize}
Let $\pi'=\pi\circ \mu\colon Y\to S$ be the induced 
proper morphism. 
Define $c$ to be the largest real number such that the log
pair $$(Y,B_Y+cF+\sum_\alpha \epsilon_\alpha G_\alpha)$$
has log canonical singularities above $X\setminus \LCS(X,B)$.
It follows from above that $F$ is supported over
$X\setminus \LCS(X,B)$, hence
$$
c=\min_\alpha \frac{1-\mult_{G_\alpha}(B_Y)-\epsilon_\alpha}
{\mult_{G_\alpha}(F)}.
$$
By perturbing the coefficients $\epsilon_\alpha$ if necessary,
we may assume that the minimum is attained at a unique prime
component of $F$, denoted by $E$. By definition, the $\R$-divisor
$$
B'=B_Y+cF+\sum_\alpha \epsilon_\alpha G_\alpha-E.
$$
does not contain $E$ in its support, and $\lceil -B'\rceil$
is effective above $X\setminus \LCS(X,B)$. The identity 
$$
(r+cm)\mu^*D-(K_Y+E+B')=r\mu^*D-K_Y-B_Y-\sum_\alpha 
\epsilon_\alpha G_\alpha+cm\bD_{m,Y}
$$
implies that the $\R$-divisor 
$
t\mu^*D-(K_Y+E+B')
$
is ample/$S$ for every $t\ge r+cm$.

Let $\mu^*D=\sum_i d_i D_i$ be the decomposition of
$\mu^*D$ into prime components. There exists a 
$\pi'$-very ample divisor $L$ on $Y$ such that 
$L+D_i$ is $\pi'$-very ample for every $i$. 
Choose {\em general} elements
$D^1_i\in \vert L+D_i\vert$ and $D^2_i\in \vert L\vert$
and set $D'=\sum_i d_i(D^1_i-D^2_i)$. The following
properties hold:
\begin{itemize}
\item[(4)] $\mu^*D\sim_\R D'$.
\item[(5)] $\langle t\mu^*D\rangle\sim \langle tD'\rangle$
for every $t\in \R$ with $\Vert \langle t\mu^*D\rangle-
t\mu^*D\Vert<1/2$.
\item[(6)] $\Supp(D')\cup \Supp(E)\cup \Supp(B')$ is included
in a simple normal crossings divisor on $Y$.
\end{itemize}
There exists $\epsilon_1>0$ such that the $\R$-divisor
$$
\langle tD'\rangle -K_Y-E-B'
$$
is ample/$S$ for $t\ge r+cm$ and 
$\Vert \langle tD'\rangle-tD'\Vert<\epsilon_1$. 
We infer by Kawamata-Viehweg vanishing (\cite{KMM}, Theorem 1-2-3) 
that the restriction map
$$
\pi'_*\cO_Y(\lceil -B'\rceil+\langle tD'\rangle)\to
\pi'_*\cO_E(\lceil -B'\vert_E\rceil+\langle tD'\vert_E\rangle)
$$
is surjective for $t\ge r+cm$ and $\Vert \langle tD'\rangle-
tD'\Vert<\epsilon_1$. On the other hand, the restrictions of
$(E,B'\vert_E)$ and $D'\vert_E$ to the generic fiber of
$E\to \pi'(E)$ satisfy the assumptions of Theorem~\ref{Knv}.
Indeed, (i) and (iii) are clear, whereas (ii) follows from the
$\pi'$-ampleness of the $\R$-divisor $(r+cm)D'\vert_E-(K_E+B'\vert_E)$.
Therefore there exists $t_0\ge r+cm$ and $0<\epsilon_2<\epsilon_1$
such that
$$
\pi'_*\cO_E(\lceil -B'\vert_E\rceil+\langle tD'\vert_E\rangle)\ne 0
$$
for every $t$ with $t\ge t_0$ and $\Vert \langle tD'\rangle-
tD'\Vert<\epsilon_2$. In particular,
$$
E\nsubseteq \Bsl_{\pi'}\vert \lceil -B'\rceil+\langle tD'\rangle
\vert \mbox { for } t\ge t_0, \Vert \langle tD'\rangle-
tD'\Vert<\epsilon_2. 
$$
Choose $0<\epsilon_3<\min(\epsilon_2,\frac{1}{2})$. We obtain
$$
E\nsubseteq \Bsl_{\pi'}\vert \lceil -B'\rceil+\langle t\mu^*D\rangle
\vert \mbox { for } t\ge t_0, \Vert \langle t\mu^*D\rangle-
t\mu^*D\Vert<\epsilon_3.
$$
It is easy to see that there exists $\epsilon_4<\epsilon_3$
such that
$$
\lceil -B'\rceil+\langle t\mu^*D\rangle\le \lceil -B'+t\mu^*D
\rceil
\mbox{ for } 
\Vert \langle t\mu^*D\rangle-t\mu^*D\Vert<\epsilon_4.
$$
By Lemma~\ref{app}, there exists a multiple $j$ of $m$ such 
that $j\ge t_0$, 
$\Vert j\mu^*D\-\langle j\mu^*D\rangle\Vert<\epsilon_4$
and $\mult_E(j\mu^*D-\langle j\mu^*D\rangle)\le 0$. In
particular, there exists a rational function $a\in 
k(Y)^\times$ such that the divisor
$$
(a)+\lceil -B'\rceil+\langle j\mu^*D\rangle
$$
is effective and it does not contain $E$ in its support.
Since $\mult_E(B')=0$, the latter property is equivalent
to 
$$
\mult_E((a)+\langle j\mu^*D\rangle)=0.
$$
On the other hand, the following inclusions hold
$$
\pi'_*\cO_Y(\lceil -B'\rceil+\langle j\mu^*D\rangle)
\subseteq \pi'_*\cO_Y(\lceil-B_Y+j\mu^*D\rceil)\subseteq
\pi'_*\cO_Y(j\mu^*D).
$$
Therefore $(a)+j\mu^*D\ge 0$. In particular,
$\mult_E(j\mu^*D-\langle j\mu^*D\rangle)\ge 0$.
The opposite inclusion holds by assumption, hence
$\mult_E(j\mu^*D)\in \Z$. Therefore the effective
$\R$-divisor $(a)+j\mu^*D$ has multiplicity zero
at $E$. In particular, $\mult_E(\bD_j-\overline{D})=0$.
By Lemma~\ref{neg}, there exists an open set $V$ in
$X$ such that $V\cap \mu(E)\ne \emptyset$ and
$(\overline{D}-\bD_j)\vert_V=0$.

{\em Step 4}. Let $m\ge 1$. By Lemma~\ref{neg}, there
exists a closed subset $\Bsl_\pi\vert mD\vert$ of
$X$ such that for every prime b-divisor $E$ of $X$,
$\mult_E(\overline{D}-\bD_j)=0$ if and only if
$c_X(E)\subset \Bsl_\pi\vert mD\vert$. For $m\vert m'$,
we have $\bD_m\le \bD_{m'}$, hence
$\Bsl_\pi\vert mD\vert \supseteq \Bsl_\pi\vert m'D\vert$.
From above, the following properties hold:
\begin{itemize}
\item[(1)] $\Bsl_\pi\vert mD\vert\subsetneq X$ for $m\ge 1$.
\item[(2)] For every $m$, there exists $m\vert m'$ such that
$\Bsl_\pi\vert mD\vert \supsetneq \Bsl_\pi\vert m'D\vert$.
\end{itemize}
By Noetherian induction, there exists a positive integer
$m$ such that $\Bsl_\pi\vert mD\vert=\emptyset$, that is 
$\bD_m=\overline{D}$. Therefore $D$ is rational, and the
Cartier divisor $mD$ is $\pi$-free.
\end{proof}

We can use Theorem~\ref{fr} to reprove Kawamata's
result (\cite{lcd}, Theorem 1) on the finite generation
of the log canonical ring modulo the existence of a 
Cutkosky-Kawamata-Moriwaki decomposition:

\begin{thm}\cite{lcd} Let $(X,B)$ be a log variety with Kawamata 
log terminal singularities such that $K+B$ is $\Q$-Cartier 
and let $\pi\colon X\to S$ be a proper morphism. Assume 
that $\bP$ is a b-nef/$S$ and b-big/$S$ $\R$-b-divisor 
of $X$ such that 
\begin{itemize}
\item[(1)] $\bP\le \overline{K+B}$.
\item[(2)] $\cR_{X/S}(\bP)=\cR_{X/S}(K+B)$.
\end{itemize}
Then $\bP$ is rational and b-semiample/$S$. In particular,
the log canonical $\cO_S$-algebra $\cR_{X/S}(K+B)$ is
finitely generated.
\end{thm}

\begin{proof} {\em Step 1}. Let $\mu:Y\to X$ be a resolution of 
singularities such that $\bP=\overline{\bP_Y}$. Let 
$\mu^*(K+B)=K_Y+B_Y$ be the log pullback, and let $B_Y=B^+_Y-B^-_Y$ 
be the decomposition of $B_Y$ into its positive and negative part, 
respectively. The effective $\Q$-divisor $B^-_Y$ is $\mu$-exceptional,
and
$$
K_Y+B^+_Y=\mu^*(K+B)+B_Y^-.
$$
In particular, $\bP_Y\le K_Y+B^+_Y$ and
$
\cR_{Y/S}(K_Y+B^+_Y)=\cR_{X/S}(K+B).
$
Therefore $(Y,B_Y^+)$ is a log variety with Kawamata log terminal
singularities, $\bP\le \overline{K_Y+B^+_Y}$, 
$\cR_{Y/S}(\bP)=\cR_{Y/S}(K_Y+B^+_Y)$. Furthermore, $\bP=\overline{\bP_Y}$.

{\em Step 2}. By Step 1, we may assume that $X$ is nonsingular,
$\bP=\overline{P}$ for some $\R$-divisor $P$ on $X$, and if
we denote $N=K+B-P$, then $\Supp(B)\cup\Supp(N)\cup\Supp(P)$
is included in a simple normal crossings divisor on $X$.
Define $B'=B-N$ and let $I$ be the smallest positive integer such 
that $I(K+B)$ is an integral divisor. For $I|j$, the following 
inequalities hold:
$$
\lceil jP-B' \rceil\le \lceil jP+N \rceil\le \lceil 
jP+jN \rceil= j(K+B).
$$
Since $\pi_*\cO_X(jP)=\pi_*\cO_X(j(K+B))$, we obtain
$$
\pi_*\cO_X(\lceil jP-B'\rceil)
\subseteq\pi_*\cO_X(jP) \mbox{ for } I|j.
$$ 
Furthermore, $P=K+B'$, hence $2P-(K+B')=P$ is nef/$S$ 
and big/$S$. By Theorem~\ref{fr}, $P$ is rational and 
semiample/$S$.
\end{proof}

%%%%%%%%%%%%%%%%%%%%%%%%%%%%%%%%%%%%%%%
%%%%%%%%%%%%%%%%%%%%%%%%%%%%%%%%%%%%%%%

\section{Finite generation of some algebras}

%%%%%%%%%%%%%%%%%%%%%%%%%%%%%%%%%%%%%%%
%%%%%%%%%%%%%%%%%%%%%%%%%%%%%%%%%%%%%%%

\begin{prop}\label{nb} Let $\pi\colon X\to S$ be a 
proper surjective morphism and let $D$ be a nef/$S$
and big/$S$ $\R$-divisor on $X$. Then: 
\begin{itemize}
\item[(1)] The characteristic limit of the $\cO_S$-algebra
$\cR_{X/S}(D)$ is $\overline{D}$.
\item[(2)] The $\cO_S$-algebra $\cR_{X/S}(D)$ is finitely 
generated if and only if $D$ is rational and semiample/$S$.
\item[(3)] Let $(X,B)$ be a log pair structure on
$X$, and let $\bD_\bullet$ be the characteristic sequence 
of $\cR_{X/S}(D)$. Then for every $j\ge 1$, there exists a 
positive integer $I(j)$ such that the natural inclusion
of sheaves
$$
\cO_X(\lceil \bA(X,B)+j\bD_i\rceil)\subseteq
\cO_X(\lceil \bA(X,B)+j\overline{D}\rceil)
$$
is an equality for $I(j)\vert i$.
\end{itemize}
\end{prop}

\begin{proof} If $\mu\colon Y\to X$ is a birational
contraction, it is clear that we may replace $X$ and
$D$ by $Y$ and $\mu^*D$, respectively. In particular,
we may assume that $X$ is nonsingular and $\pi$ is
projective. We may also assume that $S$ is affine.

We claim that there exists a sequence of ample/$S$
$\Q$-divisors $Q_k$, supported by some fixed reduced divisor 
on $X$, such that $Q_k\le D$ and 
$\lim_{k\to \infty} Q_k=D$.
Indeed, by Kodaira's Lemma, there exists an effective
$\R$-divisor $E$ on $X$ such that $D-E$ is an ample/$S$
$\Q$-divisor. Since $D$ is nef/$S$, it follows that 
$D-\frac{1}{k}E$ is an ample/$S$ $\R$-divisor for every 
$k\ge 1$. Since ampleness is an open condition and $X$ 
is $\Q$-factorial, we may find for each $k$ an ample/$S$ 
$\Q$-divisor $D_k$, supported by $\Supp(D)\cup \Supp(E)$, 
such that $Q_k\le D-\frac{1}{k}E$ and 
$\Vert D-\frac{1}{k}E-Q_k\Vert\le \frac{1}{k}$.

(1) Consider the sequence $(Q_k)_{k\ge 1}$ constructed
above. For each $k\ge 1$, there exists a positive integer 
$m_k$ such that $m_kQ_k$ is free/$S$. Since $m_kQ_k\le 
m_kD$, we obtain $\overline{m_kQ_k}\le \bM_{m_k}$, that is 
$\overline{Q_k}\le \bD_{m_k}$. In particular, 
$$
\overline{Q_k}\le \bD.
$$
Since $\lim_{k\to \infty} Q_k=D$, we obtain
$$
\overline{D}=\lim_{k\to \infty} \overline{Q_k} \le \bD.
$$
The opposite inequality always holds, hence
$
\bD=\overline{D}.
$

(2) If $\cR_{X/S}(D)$ is finitely generated, then the 
characteristic limit $\bD$ is rational and semiample/$S$. 
By (1), $D$ is rational and semiample/$S$. The converse is clear.

(3) If $\mu\colon Y\to X$ is a birational contraction,
let $\mu^*(K+B)=K_Y+B_Y$. Then $\bA(X,B)=\bA(Y,B_Y)$ and
we may replace $(X,B)$ by $(Y,B_Y)$. In particular, using
~\cite{KMM}, Corollary 0-3-6, we may assume that $X$ is 
nonsingular and contains a simple normal crossings divisor 
supporting $B,D$ and $Q_k$, for every $k\ge 1$. 

Fix $j\ge 1$. 
Since $jQ_k\le jD$ and $\lim_{k\to \infty} jQ_k=jD$,
there exists $k=k(j)$ such that
$$
\lceil -B+jQ_k\rceil=\lceil -B+jD\rceil.
$$
We have $\overline{Q_k}\le \bD_{m_k}$, hence for $m_k\vert i$ 
the following inclusion holds:
$$
\cO_X(\lceil \bA(X,B)+j\overline{Q_k}\rceil) \subseteq
\cO_X(\lceil \bA(X,B)+j\bD_i \rceil)
\subseteq \cO_X(\lceil \bA(X,B)+j\overline{D}\rceil)
$$
On the other hand, by the simple normal crossings 
assumption, we have 
$$
\cO_X(\lceil \bA(X,B)+j\overline{Q_k}\rceil)=
\cO_X(\lceil -B+jQ_k\rceil)
$$
and 
$$
\cO_X(\lceil \bA(X,B)+j\overline{D}\rceil)=
\cO_X(\lceil -B+jD\rceil).
$$
The claim holds if we set $I(j)=m_{k(j)}$.
\end{proof}
	
\begin{thm}\label{divnefbig} 
Let $(X,B)$ be a log pair with Kawamata log terminal singularities
and let $\pi\colon X\to S$ be a proper morphism. Assume that $D$ is
a nef/$S$ and big/$S$ $\R$-divisor on $X$ such that 
\begin{itemize}
\item[(1)] $\cR_{X/S}(D)$ is $\bA(X,B)$-asymptotically 
saturated.
\item[(2)] $rD-(K+B)$ is nef/$S$, for some $r\in \Q$.
\end{itemize}
Then $\cR_{X/S}(D)$ is a finitely generated $\cO_S$-algebra.
\end{thm}

\begin{proof} By Proposition~\ref{nb}, (1) is 
equivalent to the following property: 
$$
\pi_*\cO_X(\lceil \bA(X,B)+j\overline{D}\rceil)
\subseteq 
\pi_*\cO_X(jD)\mbox{ for } I\vert j.
$$ 
By Theorem~\ref{fr}, $D$ is rational and semiample/$S$. 
Therefore the algebra $\cR_{X/S}(D)$ is finitely generated.
\end{proof}

\begin{thm}\label{bX} 
Let $(X,B)$ be a log pair and let $\pi\colon X\to S$ be a 
contraction. Let $\cL$ be a bounded functional algebra 
of $X/S$, with characteristic sequence $(\bD_i)_{i\in \N(\cL)}$ 
and characteristic limit $\bD$. Assume that the following 
properties hold:
\begin{itemize}
\item[(1)] $\cL$ is asymptotically $\bA(X,B)$-saturated.
\item[(2)] $\bD_i=\overline{\bD_{i,X}}$ for every $i\in \N(\cL)$.
\item[(3)] $r\bD_X-(K+B)$ is nef/$S$ and big/$S$, for some 
$r\in \Q$.
\item[(4)] There exists $i\in \N(\cL)$ and an open neighborhood 
$U$ of $\LCS(X,B)$ such that $\bD_i\vert_U=\bD\vert_U$.
\end{itemize}
Then $\cL$ is a finitely generated $\cO_S$-algebra.
\end{thm}

\begin{proof} If $\mu\colon Y\to X$ is a birational 
contraction, let $\mu^*(K+B)=K_Y+B_Y$ be the induced
crepant log pair structure on $Y$. We have 
$\bA(X,B)=\bA(Y,B_Y)$. By (2) and (3), we have 
$\bD=\overline{\bD_X}$. Therefore the asumptions are
preserved if we replace $(X,B)$ by $(Y,B_Y)$.
In particular, we may assume that $X$ is nonsingular and 
admits a normal crossings divisor supporting $B$ and 
every $\bD_{i,X}$. 

By (2), we have $\bD_i=\overline{\bD_{i,X}}$ for every $i$.
In particular, every $\Q$-divisor $\bD_{i,X}$ is nef/$S$ 
and their limit $D=\lim_{i\to \infty} \bD_{i,X}$ is a nef/$S$ 
$\R$-divisor.
By the normal crossings assumption and (2), asymptotic 
saturation is equivalent to the following inclusions
$$
\pi_*\cO_X(\lceil -B+j\bD_{i,X}\rceil)\subseteq \pi_*\cO_X(\bM_j)
\mbox{ for } I\vert i,j.
$$
For each $j$, we have $\lceil -B+j\bD_{i,X}\rceil=
\lceil -B+jD\rceil$ for $i$ sufficiently large and 
divisible. Therefore
$$
\pi_*\cO_X(\lceil -B+jD\rceil)
\subseteq \pi_*\cO_X(\bM_j) \subseteq \pi_*\cO_X(jD) 
\mbox{ for } I|j.
$$
Since $D-\bD_{i,X}$ is effective and zero on the open 
neighborhood $U$ of $\LCS(X,B)$, and $i\bD_{i,X}$
is free/$S$, the $\R$-divisor $iD$ is free/$S$ on $U$. 
Furthermore, $rD-(K+B)$ is nef/$S$ and big/$S$. 
Theorem~\ref{fr} implies that $D$ is a semiample/$S$ $\Q$-divisor. 

Let $j$ be a multiple of $I$ such that $jD$ is free/$S$. 
There exists a contraction $\Phi\colon X\to T$, defined over
$S$, and an ample/$S$ Cartier divisor $A$ on $T$ such that 
$jD\sim\Phi^*(A)$. 

Set $C=\lceil -B^{\ge 1}\rceil\le 0$, so that $\LCS(X,B)=\Supp(C)$. 
By assumption, $\bD_j=\bD$ over the neighborhood $U$ of 
$\Supp(C)$. By Lemma~\ref{neg'}, there exist an open 
set $V_j\supset \Phi(\Supp(C))$ such that 
$$
\bD_j\vert_{\Phi^{-1}(V_j)} =\bD\vert_{\Phi^{-1}(V_j)}.
$$
If $\Phi(\Supp(C))=T$, then $V_j=T$ and we are done. Assume now 
that $\Phi(\Supp(C))\subsetneqq T$.

Since $A$ is ample/$S$, we may replace $j$ by a multiple, so 
that the subsheaf
$$
\Phi_*\cO_Y(C)\otimes \cO_T(jA)\subset \cO_T(jA)
$$
is relatively generated by global sections. In particular,
the divisor $\lceil C+jD\rceil$ is free on 
$Y\setminus \Phi^{-1}\Phi(\Supp(C))$. By asymptotic 
saturation, we have
$$
\pi_*\cO_X(\lceil C+jD\rceil)\subseteq \pi_*\cO_X(\lceil
-B+jD\rceil)\subseteq \pi_*\cO_X(\bM_j). 
$$
Therefore $\overline{jD}\le j\bD_j$ above 
$X\setminus \Phi^{-1}\Phi(\Supp(C))$. The opposite
inclusion always holds, hence 
$\bD_j=\bD$ over $X \setminus \Phi^{-1}\Phi(\Supp(C))$.
The two open sets $\Phi^{-1}(V_j)$ and 
$X \setminus \Phi^{-1}\Phi(\Supp(C))$ cover $X$, hence
we obtain $\bD_j=\bD$. 

Therefore the $\cO_S$-algebra $\cL$ is finitely generated.
\end{proof}

%%%%%%%%%%%%%%%%%%%%%%%%%%%%%%%%%%%%%%%
%%%%%%%%%%%%%%%%%%%%%%%%%%%%%%%%%%%%%%%

\end{document}